\pdfoutput=1
\documentclass{amsart}
\usepackage{amsmath}
\usepackage{amsthm}
\usepackage{microtype}
\usepackage[hmargin=4cm]{geometry}
\usepackage{hyperref}
\newtheorem{theorem}{Theorem}
\newtheorem{lemma}{Lemma}

\begin{document}

\title{Generating functions for plateaus in Motzkin paths}

\author{Dan Drake}
\address{Department of Mathematical Sciences\\
  Korea Advanced Institute of Science and Technology\\
  Daejeon, Korea}
\email{ddrake@member.ams.org}
\urladdr{http://mathsci.kaist.ac.kr/~drake}

\author{Ryan Gantner}
\address{Department of Mathematical and Computing Sciences\\
  St.\ John Fisher College\\
  Rochester, New York}
\email{rgantner@sjfc.edu}
\urladdr{http://citadel.sjfc.edu/faculty/rgantner/}


\begin{abstract}
  A plateau in a Motzkin path is a sequence of three steps: an up step,
  a horizontal step, then a down step. We find three different forms for
  the bivariate generating function for plateaus in Motzkin paths, then
  generalize to longer plateaus. We conclude by describing a further
  generalization: a continued fraction form from which one can easily
  derive new multivariate generating functions for various kinds of path
  statistics. Several examples of generating functions are given using
  this technique.
\end{abstract}

\maketitle

A \emph{Motzkin path} is a finite sequence of steps with the following
properties: each step is ``up'' (labeled U), ``horizontal'' (labeled
H), or ``down'' (labeled D); at any point in the sequence, the
number of up steps is at least as big as the number of down steps; the
total number of up steps in the sequence is equal to the total number of
down steps. The number of steps in the sequence is referred to as the
length of the path. Motzkin paths can be visualized as graphs in the
plane beginning at $(0,0)$ and consisting of up steps in direction
$(1,1)$, horizontal steps in direction $(1,0)$, and down steps in
direction $(1,-1$). In this visualization, the path ends at point
$(n,0)$, where $n$ is the length of the path, and never passes below the
horizontal axis. We will use this visualization to refer to the points
between steps as vertices.

To begin, we define a \emph{plateau} in a Motzkin path to be a
subsequence of the path consisting of an up step immediately followed by
a horizontal step immediately followed by a down step (such a
subsequence is often abbreviated UHD). In this paper we establish several forms for generating functions that count plateaus for Motzkin
paths. We then generalize those methods to longer plateaus and also use
continued fractions to derive similar bi- and multivariate generating
functions.

The enumeration of plateaus in Motzkin paths has applications in the
theory of RNA secondary structures \cite{nebel, viennot}. Work by
Prodinger and Wagner \cite{prodinger} examines plateaus in Motzkin paths
with a different goal.

\section{Forms for generating functions}

In this section we state three theorems about forms for the generating function which counts plateaus for Motzkin paths.
Let $c_n^p$ be the number of Motzkin paths of length $n$ with $p$
plateaus, and let
\[
g(x,y) = \sum_{n=0}^{\infty} \sum_{p=0}^{\lfloor \frac{n}{3} \rfloor} c_n^p x^n y^p ,
\]
so $g$ is the bivariate generating function for this array.

\begin{theorem}
\label{approach1theorem}
The function $g$ has an integral/differential form given by
\begin{equation}
\label{approach1result}
g(x,y) = \frac{1-2x^{3}}{1-2x^{3}(1-y)} \left(
  f_{0}(x) + \frac{x}{1-2x^{3}} \frac{\partial}{\partial x} x^{3} \int_0^y
  g(x,t)\,\mathrm{d}t \right),
\end{equation}
where
\[
f_{0}(x) =  \frac{1 - x + x^3 - \sqrt{(1-
    x + x^3)^2 - 4x^2}}{2x^2}.
    \]
\end{theorem}

\begin{theorem}
\label{approach2theorem}
The function $g$ has a differential form given by
\begin{equation}
\label{approach2result2}
\frac{\partial}{\partial x} x g\left(x, \frac{z}{x^3}\right)
= (1-z-2x^3)\frac{\partial}{\partial z} g\left(x, \frac{z}{x^3}\right).
\end{equation}
\end{theorem}

\begin{theorem}
\label{approach3theorem}
The function $g$ has an explicit form given by
\begin{equation}
\label{gfexplicit}
g(x,y) = \frac{1 - x + x^3 - x^3y - \sqrt{(1 - 3x + x^3 - x^3 y)(1 + x +
    x^3 - x^3 y)}}{2x^2}.
\end{equation}
\end{theorem}

\section{Proof of Theorems \ref{approach1theorem}--\ref{approach3theorem}}

Our plateau-counting formulas depend on a recursion among the $c_n^p$.
This can be obtained by the process of ``sewing in'' a plateau, as the
first lemma reveals. After stating and proving the lemma, we'll turn to
proving the theorems.

\begin{lemma}
\label{recursionlemma}
If we set $c_n^p = 0$ when $p$ is negative or $n$ is negative, then the $c_n^p$ satisfy
\begin{equation}
    \label{recursion}
c_n^p = \frac{n-2p}{p} c_{n-3}^{p-1} + 2 c_{n-3}^p
\end{equation}
for all
$n$ and for all $p > 0$.
\end{lemma}

\begin{proof}[Proof of Lemma \ref{recursionlemma}]
To get a path of length $n$ with $p$ plateaus, we can start
with a path of length $n-3$ with $p-1$ plateaus and sew in a plateau by
inserting a UHD subpath at any vertex. There are $n-3+1$ vertices into
which a plateau can be sewn in a path of length $n-3$. Sewing in a
plateau in this way always gives us a new plateau and a Motzkin path of
length $n$. Sometimes, however, the sewing operation destroys an
existing plateau in the path of length $n-3$. If a plateau is sewn on a
vertex that is adjacent to the horizontal step in an existing plateau
(such a vertex will be said to be \emph{inside} a plateau), we destroy
the plateau that was in the original path of length $n-3$ in order to
create one in the path of length $n$. In all other places, sewing in a
plateau creates a new plateau without destroying a previous one. In
summary, to get a path of length $n$ with $p$ plateaus, we can sew in a
plateau at any of the $n-3+1-2(p-1)$ vertices which are not inside a
plateau in a path of length $n-3$ with $p-1$ plateaus, or we can sew in
a plateau at any of the $2p$ vertices inside plateaus of a path of
length $n-3$ with $p$ plateaus. By doing this, we can generate all of
the paths of length $n$ with $p$ plateaus. In fact, we generate each
path $p$ times, according to our choice of which plateau is sewn in.
Thus, we get the recursion
\[
    c_n^p = \frac{n-2-2(p-1)}{p} c_{n-3}^{p-1} +
     \frac{2p}{p}c_{n-3}^p
\]
which simplifies to~\eqref{recursion}.
\end{proof}

As a note, to use the recursion, we'll need to be given the sequence $\{c_{n}^{0}\}$. We'll see that the sequence $\{c_{n}^{0}\}$ obeys a
Catalan-like recurrence relation given below in \eqref{noplats}.

To give a sense of what this array of numbers looks like, see the table
of values of $c_n^p$ for small $n$ and $p$ given in
Table~\ref{smallvalues}. That triangle is sequence A114583 in the
OEIS~\cite{oeis}.

\begin{table}
\centering
\begin{tabular}{r|rrrrr}
     & $p$: 0 & 1 & 2 & 3 & 4 \\
\hline
$n$: 0 & 1 &   &   &   &   \\
1    & 1 &   &   &   &   \\
2    & 2 &   &   &   &   \\
3    & 3 & 1 &   &   &   \\
4    & 7 & 2 &   &   &   \\
5    & 15& 6 &   &   &   \\
6    &36 &14 & 1 &   &   \\
7    &85 &39 & 3 &   &   \\
8    &209&102& 12&   &   \\
9    &517&280& 37& 1 &   \\
10   &1303& 758& 123& 4  \\
11   &3312&2085& 381& 20 \\
12   &8510& 5730& 1194& 76& 1 \\
13   &22029& 15849& 3657& 295& 5\\
14   &57447& 43914& 11187& 1056& 30
\end{tabular}
\caption{
$c_n^p$ for small values of $n$ and $p$.
}\label{smallvalues}
\end{table}


\subsection{Proof of Theorem \ref{approach1theorem}: generating functions for each
  \texorpdfstring{$p$}{\it p}}
\label{section:approach1}

To prove Theorem~\ref{approach1theorem}, we first find the
generating function for the sequence $\{ c_n^0\}_{n=0}^{\infty}$, the
first column of Table~\ref{smallvalues}. Let $f_{0}(x)$ denote that
generating function. We first show the following.

\begin{lemma}
\label{columnlemma}
The function $f_{0}$ is given by
\begin{equation}
\label{genfunD}
f_{0}(x) = \sum_{n=0}^{\infty} c_n^0 x^n = \frac{1 - x + x^3 - \sqrt{(1-
    x + x^3)^2 - 4x^2}}{2x^2}.
\end{equation}
Furthermore, the sequence $\{c_{n}^{0}\}$ satisfies
\begin{equation}
\label{noplats}
     c_n^0 = c_{n-1}^0 + c_{n-2}^0 + \sum_{k=2}^{n-2} c_k^0 c_{n-k-2}^0.
\end{equation}
\end{lemma}

\begin{proof}[Proof of Lemma \ref{columnlemma}]
The terms of the sequence $\{c_{n}^{0}\}$ are the number of
Motzkin paths of various lengths which have no plateaus. Each
such path falls into one of three categories.

\emph{Category 0: the empty path.} This contributes $1$ to $f_{0}(x)$.

\emph{Category 1: the path starts with a horizontal step.} In this case,
the horizontal step is followed by a Motzkin path of length $n-1$ with
no plateaus. The generating function for such paths is $x f_{0}(x)$.

\emph{Category 2: the path starts with an up step.} In this case, the
path can be decomposed into
\begin{enumerate}
\item an up step,
\item a plateau-free Motzkin path that is not a single horizontal step,
\item a down step, then
\item any plateau-free Motzkin path.
\end{enumerate}
That decomposition tells us that the generating function for paths in
category 2 is $x (f_{0}(x) - x) x f_{0}(x)$.

Every plateau-free Motzkin path falls into exactly one of the above
categories, so
\[
f_{0}(x) = 1 + x f_{0}(x) + x^{2} f_{0}(x) (f_{0}(x) - x),
\]
which after an application of the quadratic formula yields
\[
    f_{0}(x) = \frac{1-x + x^3 \pm \sqrt{(1-x + x^3)^2 - 4x^2}}{2x^2}.
\]
The limit of $f_{0}(x)$ as $x$ goes to zero is $1$ (the empty path has
no plateau), so we may take limits of both possibilities for $f_{0}(x)$
to see that only the subtraction term in the numerator makes sense, which gives~\eqref{genfunD}.

Similar reasoning with the above categories and decomposition yields the
recurrence relation for $c_{n}^{0}$ in~\eqref{noplats}.
\end{proof}

\emph{Remark}.  Note that the function in~\eqref{genfunD} is the generating function for OEIS sequence A114584~\cite{oeis}.

\begin{proof}[Proof of Theorem \ref{approach1theorem}]
  Knowing $f_{0}(x)$ and the recurrence relation from
  Lemma~\ref{recursionlemma}, we can calculate the generating functions
  $f_{p}(x) = \sum_n c_{n}^{p} x^{n}$. For $p > 0$, one can multiply
  both sides of~\eqref{recursion} by $x^{n}$ and sum over $n \ge 0$ to
  get
\begin{align*}
  f_{p}(x) = \sum_{n=0}^{\infty} c_{n}^{p} x^{n}
  &= \sum_{n=0}^{\infty} \frac{n-2p}{p} c_{n-3}^{p-1} x^{n} + 2
  \sum_{n=0}^{\infty} c_{n-3}^{p} x^{n}\\
  &= \frac{1}{p} \sum_{n=0}^{\infty} n c_{n-3}^{p-1} x^{n}
     - 2 x^{3} f_{p-1}(x) + 2 x^{3} f_{p}(x).
\end{align*}
The remaining sum above equals
\[
\frac{x}{p} \frac{\mathrm{d}}{\mathrm{d}x} \left[ x^{3} f_{p-1}(x)\right],
\]
so one solves to find a recurrence relation for $f_{p}$ for $p \ge 1$:
\begin{equation}
  \label{fp-gf}
  f_{p}(x) = \frac{1}{1-2x^{3}} \left( \frac{x}{p}
    \frac{\mathrm{d}}{\mathrm{d}x} \left[x^{3} f_{p-1}(x)\right] - 2 x^{3} f_{p-1}(x) \right),
\end{equation}
which can also be written
\begin{equation}
  \label{fp-gf2}
  f_{p}(x) = \frac{x^{2p+1}}{p(1-2x^{3})} \frac{\mathrm{d}}{\mathrm{d}x} \left[x^{3-2p} f_{p-1}(x)\right].
  \tag{\ref*{fp-gf}$'$}
\end{equation}

We can use the recurrence for the generating functions $f_{p}(x)$ to
derive the integral/differential form for $g(x,y)$ in~\eqref{approach1result}. Start
with~\eqref{fp-gf}, multiply both sides by $y^{p}$ and sum over all $p
\ge 1$. This yields
\[
g(x, y) - f_{0}(x) = \frac{1}{1-2x^{3}} \sum_{p=1}^{\infty} \left(
  \frac{x}{p} \frac{\partial}{\partial x} x^{3} f_{p-1}(x) - 2x^{3}
  f_{p-1}(x) \right) y^{p}.
\]
A little rearrangement and simplification leads to
\[
g(x, y) = f_{0}(x) + \frac{y}{1-2x^{3}} \left(
          \sum_{p=1}^{\infty} \frac{x}{p} \frac{\partial}{\partial x}
          x^{3} f_{p-1}(x) y^{p-1} - 2x^{3} g(x,y) \right),
\]
and since
\[
\sum_{p=1}^{\infty} \frac{1}{p} \frac{\partial}{\partial x}
   x^{3} f_{p-1}(x) y^{p-1}
=
\frac{1}{y}\frac{\partial}{\partial x} x^{3} \int g(x,y)\,\mathrm{d}y,
\]
we can solve for $g$ to get our first expression for $g(x,y)$ to get~\eqref{approach1result}.
\end{proof}

\subsection{Proof of Theorem~\ref{approach2theorem}: summing array diagonals}

In Theorem~\ref{approach1theorem}, we found a functional equation for $g$ by finding
generating functions for the columns of Table~\ref{smallvalues}. For Theorem~\ref{approach2theorem},
we find a form for $g$ by looking at the diagonals of that table. Let
$h_k(z) = \sum_{m=0}^{\infty} d_m z^m$, where $d_m = c_{3m+k}^m$. For
instance, $h_0(z) = \sum_{m=0}^{\infty} c_{3m}^m z^m =
\sum_{m=0}^{\infty} 1 z^m = 1/(1-z)$, since there is just one path of
length $3m$ with $m$ plateaus. Using the recursion from Lemma~\ref{recursionlemma}, we
get
\[
c_{3m+k}^m = \frac{m+k}{m} c_{3m+k-3}^{m-1} + 2 c_{3m+k-3}^m.
\]
Substituting this into the generating function $h_k$, we get
\begin{align*}
\sum_{m=0}^{\infty} d_m z^m
&= c_k^0 + \sum_{m=1}^{\infty} \left( \frac{m+k}{m} c_{3(m-1)+k}^{m-1} + 2c_{3m + (k-3)}^m \right) z^m \\
&= c_k^0 + \sum_{m=1}^{\infty} c_{3(m-1)+k}^{m-1} z^m + \sum_{m=1}^{\infty} \frac{k}{m} c_{3(m-1)+k}^{m-1} z^m +
   2 \sum_{m=1}^{\infty} c_{3m + (k-3)}^m  z^m
\end{align*}
which reduces to
\[
 h_k(z) = z h_k(z) + k \int_0^z h_k(t)\, \mathrm{d}t + 2 h_{k-3}(z) -2c_{k-3}^0 + c_k^0.
\]
Upon differentiation, we get the differential difference equation
\begin{equation}
\label{difdif}
 h'_k(z) = z h'_k(z) + h_k(z) + k h_k(z) + 2 h'_{k-3}(z)
\end{equation}

Returning to $g$, we can use \eqref{difdif} to derive~\eqref{approach2result2}:
first observe that
\[
g(x,y) = \sum_{k \ge 0} h_k(x^3 y) x^k,
\]
so that
\begin{gather*}
\frac{\partial}{\partial x} x g\left(x, \frac{z}{x^3}\right)
= \sum_{k \ge 0} (k+1) h_k(z) x^k \quad\text{and}\\
\frac{\partial}{\partial z} g\left(x, \frac{z}{x^3}\right)
= \sum_{k \ge 0} h'_k(z) x^k.
\end{gather*}
Since \eqref{difdif} can be rearranged into
\[
(k+1)h_k(z) = (1-z) h'_k(z) - 2 h'_{k-3}(z),
\]
if one multiplies both sides of that equation by $x^k$ and sums over all
$k$, we obtain~\eqref{approach2result2} and complete the proof of Theorem \ref{approach2theorem}.
\qed

\emph{Remark.}  Much more can be said about the functions $h_k$ from the proof above.  For instance,
one may solve the differential equation \eqref{difdif} using standard techniques to
obtain, for $k \ge 3$,
\begin{equation}
  \label{hrecurrence}
  h_{k}(z) = \frac{1}{(1-z)^{k+1}}\left( c_{k}^{0} +
    2 \int_{0}^{z} (1-t)^{k} h_{k-3}'(t)\,\mathrm{d}t \right).
\end{equation}
We just need three initial conditions: when $k=0$, we already know that
$h_{0}(z) = 1/(1-z)$. If $k=1$, the recurrence yields the sequence of
positive natural numbers, whose generating function is $1/(1-z)^{2}$,
and if $k=2$, the recurrence relation produces the sequence $2, 6, 12,
20, 30, 42,\dots$, whose generating function is simply the derivative
of the previous: $2/(1-z)^{3}$. With these three initial generating
functions and the recurrence relation \eqref{hrecurrence}, any
$h_{k}(z)$ can be found. Note that if one defines $h_{k}(z) = 0$ for
negative $k$, the recurrence \eqref{hrecurrence} holds for all
nonnegative $k$.

Furthermore, for $k = 0$, $1$, and $2$, $h_{k}$ is a rational function whose
denominator is $(1-z)^{k+1}$, and it is a simple matter to use
induction to prove that the same is true for all $k$: if $h_{k-3}$ is
a rational function of the form $p/(1-z)^{k-2}$, where $p$ is a
polynomial of degree $d$, then the integrand in \eqref{hrecurrence} is
a polynomial of degree $d+1$, so the second factor in that expression
is a polynomial of degree $d+2$, and we see that the numerator of
$h_k$ has degree $2 \lfloor k/3 \rfloor$.  Moreover, if $N_k$ denotes
the numerator of $h_k$, the numerators obey the differential
difference equation
\[
N_{k}(z) = c_{k}^{0} + 2 \int_{0}^{z}
(1-t)\left((1-t)\frac{\mathrm{d}}{\mathrm{d}t} N_{k-3}(t) +
  (k-2)N_{k-3}(t)\right) \, \mathrm{d}t.
\]

\subsection{Proof of Theorem \ref{approach3theorem}: explicit form}

As a third approach to finding the generating function for the array
$\{c_n^p\}$, we can take the categories of paths defined in
\autoref{section:approach1} for plateau-free Motzkin paths and use them
to work on all Motzkin paths. Each Motzkin path falls into one of three
categories: the path is empty, it begins with a horizontal step, or it
begins with an up step. If it begins with an up step, it is of the form
``U$P$D$Q$'' where $P$ and $Q$ are Motzkin paths. Since $Q$ can be any
Motzkin path, the generating function describing the possibilities for
$Q$ is simply $g$. For $P$, we can have any path, but if we use a path
that is a single horizontal step, we will create an extra plateau which
is unaccounted for in the generating function multiplication. To combat
this, we can subtract $x$ from $g$ and add $xy$ to count the plateau
properly. Therefore, $g(x,y)$ satisfies
\begin{equation}
\label{gfrelation}
g(x,y) = 1 + x g(x,y) + x (g(x,y) - x + xy) x g(x, y).
\end{equation}
We can again use the quadratic formula in that functional equation and
take limits to find the explicit form given in \eqref{gfexplicit}.
\qed

\section{Generalization to longer plateaus}

We can generalize the approaches in the three theorems above to longer plateaus. Define
a \emph{plateau of length} $r$ to be a subsequence of a Motzkin path
consisting of an up step, immediately followed by $r$ consecutive
horizontal steps, then a down step (we will abbreviate such a
subsequence as ``$\text{UH}^r \text{D}$''). Using this definition, the
previous calculations have been for plateaus of length 1. Since the groundwork has been laid, in this section we quickly state
and prove results similar to the previously stated lemmas and theorems.  We set ${_r}c_n^p$ to be the number of Motzkin paths of
length $n$ with $p$ plateaus of length $r$

\begin{lemma}
\label{newrecursionlemma}
With  ${_r}c_n^p$ as above, we have
\begin{equation}
\label{newrecursion}
{_r}c_n^p = \frac{n - (r+1)p}{p} {_r}c_{n-(r+2)}^{p-1} + (r+1) {_r}c_{n - (r+2)}^p .
\end{equation}
\end{lemma}

\begin{proof}
  We observe that there are $r+1$ vertices inside each plateau for which
  sewing in a plateau destroys the existing plateau upon creating
  another. Also, sewing in a plateau of length $r$ increases the length
  of the path by $r+2$. Therefore, we get the recursion
\[ {_r}c_n^p = \frac{n - (r+1) - (r+1)(p-1)}{p} {_r}c_{n-(r+2)}^{p-1} + \frac{(r+1)p}{p} {_r}c_{n - (r+2)}^p,\]
which simplifies to~\eqref{newrecursion}.
\end{proof}

We now quickly address each of the three approaches we used earlier for
finding forms for the generating function ${_r}g(x,y) =
\sum_{n=0}^{\infty} \sum_{p = 0}^{\lfloor n/(r+2)\rfloor} {_r}c_n^p x^n
y^p$.

\begin{theorem}
\label{approach1rtheorem}
The function  ${_r}g$ satisfies the integral/differential form
\begin{multline}
\label{approach1rresult}
{_r}g(x,y) = \frac{1-(r+1)x^{r+2}}{1-(r+1)x^{r+2}(1-y)}
  \Bigg( {_r}f_{0}(x) \\
+ \frac{x}{1-(r+1)x^{r+2}}
    \frac{\partial}{\partial x} x^{r+2} \int_0^y {_r}g(x,t)\,\mathrm{d}t
  \Bigg),
\end{multline}
where
\[
{_r}f_{0}(x) = \frac{-x + x^{2+r} + 1 - \sqrt{(x - x^{2+r} - 1)^2 - 4x^2}}{2x^2}.
\]
\end{theorem}

\begin{proof}
We find a relationship for the generating
functions ${_r}f_p(x) = \sum_{n=0}^{\infty} {_r}c_n^p x^n$ for each $p$, as we did in Lemma \ref{columnlemma}.
We use the same categorization as discussed in the proof of that lemma, simply noting that
now the second category consists of Motzkin paths that have an up step,
then a Motzkin path with no plateaus of length $r$ which is not a
sequence of $r$ horizontal steps, then a down step, then any Motzkin
path with no plateaus of length $r$. We use this to obtain
\[
{_r}f_0(x) = 1 + x {_r}f_0(x) + x^2 {_r}f_0(x) ({_r}f_0(x) - x^r),
\]
where ${_r}f_0(x)$ is the generating function for number of Motzkin
paths of length $n$ with no plateaus of length $r$. Solving the
quadratic equation above yields an explicit form
\begin{equation}
\label{genfunDr}
{_r}f_{0}(x) = \sum_{n=0}^{\infty} {_r}c_n^0 x^n = \frac{-x + x^{2+r} + 1 - \sqrt{(x - x^{2+r} - 1)^2 - 4x^2}}{2x^2}.
\end{equation}
The recursion for ${_r}c_n^p$ analogous to~\eqref{noplats} in Lemma \ref{columnlemma} is
\[
{_r}c_n^0 = {_r}c_{n-1}^0 + \sum_{k=0}^{n-2} {_r}c_k^0 {_r}c_{n-k-2}^0 - {_r}c_{n-r-2}^0.
\]
We then continue as in the proof of Theorem \ref{approach1theorem} by multiplying both sides
of~\eqref{newrecursion} by $x^n$ and summing over $n\geq 0$ to get
\[ {_r}f_{p}(x) = \sum_{n=0}^{\infty} {_r}c_{n}^{p} x^{n}
  = \sum_{n=0}^{\infty} \frac{n-(r+1)p}{p} {_r}c_{n-(r+2)}^{p-1} x^{n} + (r+1)\sum_{n=0}^{\infty} {_r}c_{n-(r+2)}^{p} x^{n}.
\]
This can be manipulated using the same techniques as earlier to get the
integral/differential form for ${_r}g$ in \eqref{approach2rresult}.
\end{proof}

Next, we can use the technique of summing array diagonals to get another form for ${_r}g$.
\begin{theorem}
\label{approach2rtheorem}
The function  ${_r}g$ satisfies
\begin{equation}
\label{approach2rresult}
\frac{\partial}{\partial x} x {_r}g\left(x, \frac{z}{x^{r+2}}\right)
= (1-z-(r+1)x^{r+2})\frac{\partial}{\partial z} {_r}g\left(x, \frac{z}{x^{r+2}}\right).
\end{equation}
\end{theorem}

\begin{proof}
Define
${_r}h_k(z)$ to be $\sum_{m=0}^{\infty} {_r}d_m z^m$, where ${_r}d_m =
c_{(r+2)m+k}^m$. Using the recursion in Lemma \ref{newrecursionlemma}, we get
\[
{_r}c_{(r+2)m+k}^m = \frac{m+k}{m} {_r}c_{(r+2)m+k-(r+2)}^{m-1} +
(r+1) {_r}c_{(r+2)m+k-(r+2)}^m,
\]
which upon substitution into the definition of ${_r}h_k$, becomes
\[ {_r}h_k(z) = {_r}c_k^0 + \sum_{m=1}^{\infty} \left( \frac{m+k}{m} {_r}c_{(r+2)(m-1)+k}^{m-1} + (r+1) {_r}c_{(r+2)m + k-(r+2)}^m \right) z^m.\]
Upon differentiation we get
\begin{equation}
\label{newdifdif}
 {_r}h'_k(z) = z {_r}h'_k(z) + {_r}h_k(z) + k {_r}h_k(z) + (r+1){_r}h'_{k-(r+2)}(z),
\end{equation}
a differential difference equation similar to~\eqref{difdif}. This gives
rise to \eqref{approach2rresult}, which is similar to~\eqref{approach2result2} from Theorem \ref{approach2theorem}.
\end{proof}

Finally, following the approach of the proof of Theorem \ref{approach3theorem}, we show the following.
\begin{theorem}
\label{approach3rtheorem}
The function ${_r}g$ has explicit form
\begin{equation}
\label{newgfexplicit}
{_r}g(x,y) = \frac{1 - x + x^{r+2} - x^{r+2} y - \sqrt{(1 - x + x^{r+2} - x^{r+2} y)^2 - 4x^2}}{2x^2}.
\end{equation}
\end{theorem}

\begin{proof}
  We understand that either a Motzkin path is empty, it begins with a
  horizontal step, or it begins with an up step. If it begins with an up
  step, it is of the form ``U$P$D$Q$'' where $P$ and $Q$ are Motzkin
  paths. Since $Q$ can be any Motzkin path, the generating function for
  that part of the decomposition is ${_r}g$. For $P$, we can have any
  path, but if we use a path that is a sequence of $r$ horizontal steps
  we will create an extra plateau which is unaccounted for in the
  generating function multiplication. To combat this, we can subtract
  $x^r$ from $g(x,y)$ and add $x^ry$. Therefore, ${_r}g(x,y)$ satisfies
\begin{equation}
\label{newgfrelation}
{_r}g(x,y) = 1 + x {_r}g(x,y) + x^{2} {_r}g(x,y) ({_r}g(x,y) - x^r + x^ry).
\end{equation}
We use the quadratic formula in~\eqref{newgfrelation} to solve for
${_r}g$ explicitly to get \eqref{newgfexplicit}.
\end{proof}

\emph{Remark}.
Observe that if $r = 0$, then the ``plateaus'' are just UD
subpaths---that is, they are peaks, and when $r=0$ the generating
function in Theorem \ref{approach3theorem} does correctly count Motzkin paths in which peaks have weight
$y$; see OEIS sequence A097860 \cite{oeis}.

\section{Further generalizations and continued fraction expansions}


\newsavebox{\mytop}
\newsavebox{\mybot}
\newlength{\toplen}
\newlength{\botlen}
\newlength{\extra}

\newcommand{\hcfrac}[2]{%
\sbox{\mytop}{$#1$}%
\sbox{\mybot}{$#2$}%
\settowidth{\toplen}{\usebox{\mytop}}%
\settowidth{\botlen}{\usebox{\mybot}}%
\setlength{\extra}{\botlen}%
\addtolength{\extra}{-\toplen}%
\setlength{\extra}{.5\extra}%
\frac{\hspace{\extra}\hspace{2.7778pt}\usebox{\mytop}\hspace{\extra}\hspace{-2.7778pt}\mid}%
{\mid\usebox{\mybot}}}

The functional equation~\eqref{newgfrelation} leads us to a continued
fraction form that generalizes the above generalization. Instead of
choosing a specific plateau length, let us rewrite~\eqref{newgfrelation}
as:
\begin{equation}
\label{newergfrelation}
G(x,y) = 1 + xG(x,y) + x (G(x,y) + C) x G(x,y);
\end{equation}
here $C$ simply stands for whatever correction we need to make for
plateaus. Rearrange that functional equation and we have:
\[
G(x, y) = \frac{1}{1 - x - x^{2} C - x^{2} G(x,y)}.
\]
Inductively replacing $G$ with the right-hand side of that equation
yields the continued fraction expansion
\begin{equation}
\label{eq:contfracexpansion}
G(x,y) = \frac{1}{1 - x - x^{2} C -
         \cfrac{x^{2}}{1 - x - x^{2} C -
         \cfrac{x^{2}}{1 - x - x^{2} C -
         \cfrac{x^{2}}{1 - \cdots\vphantom{\cfrac{x^{2}}{1}}}}}},
\end{equation}
which we will write more compactly as
\[
G(x, y) = \hcfrac{1}{1 - x - x^{2}C} - \hcfrac{x^{2}}{1 - x - x^{2}C} -
\hcfrac{x^{2}}{1 - x - x^{2}C} -\cdots
\]
We can use this discussion to get continued fraction expansions for $g$ and ${_r}g$ from earlier.

\begin{theorem}
  The function $g$ of Theorems
  \ref{approach1theorem}--\ref{approach3theorem} has a continued fraction
  expansion given by
\begin{equation}
\label{eq:cf}
g(x,y) =  \hcfrac{1}{1 - x - x^{2}(xy-x)} - \hcfrac{x^{2}}{1 - x - x^{2}(xy-x)}
 -\cdots .
\end{equation}
\end{theorem}
\begin{proof}
In our original problem, we needed to give weight $xy$ to the horizontal
step in a UHD subpath, not $x$, so in the above continued fraction, we
simply set $C = xy - x$ in \eqref{eq:contfracexpansion}.
\end{proof}

\begin{theorem}
  The function ${_r}g$ of Theorems
  \ref{approach1rtheorem}--\ref{approach3rtheorem} has a continued
  fraction expansion given by
\begin{equation}
\label{eq:cf-generalized}
{}_{r}g(x,y) =  \hcfrac{1}{1 - x - x^{2}(x^{r}y-x^{r})} - \hcfrac{x^{2}}{1 - x - x^{2}(x^{r}y-x^{r})} - \cdots
\end{equation}
\end{theorem}
\begin{proof}  Set $C = x^{r}y - x^{r}$ in \eqref{eq:contfracexpansion}.
\end{proof}

This form of the generating function has two advantages: first, we can
specialize $C$ to represent whatever plateau (or other features of the
path) that we want, and second, we now have ``infinitely many $C$s''.
The $C$ appearing at the $k$th level of~\eqref{eq:contfracexpansion}
corresponds to corrections made at height $k$ in a path (see Flajolet
\cite[Theorem 1]{flajolet:combinatorial}), so we can generalize that
continued fraction expansion to
\begin{equation}
\label{eq:contfracexpansion2}
G(x,y) = \hcfrac{1}{1 - x - x^{2}C_{1}} - \hcfrac{x^{2}}{1 - x - x^{2}C_{2}} -
\hcfrac{x^{2}}{1 - x - x^{2}C_{3}} -\cdots
\end{equation}
Now $C_{k}$ refers to the ``correction term'' for parts of the path
occurring at height $k$. We can now easily find the generating function
for many variations of the problems considered here. Here are several
examples.

\begin{theorem}
The generating function for Motzkin paths with no peaks (UD
  subpaths) is
  \[
  \frac{x^{2} - x + 1 -\sqrt{x^{4} - 2 x^{3} - x^{2} - 2 x + 1}}{2 x^{2}}.
  \]
\end{theorem}
\begin{proof} Set $C=-1$ in~\eqref{newergfrelation}
  to correct for the empty path between the up and down steps and solve
  as before.
  \end{proof}

\emph{Remark}.  This is sequence A4148 in the OEIS~\cite{oeis}.

\begin{theorem}
The generating function for Motzkin paths in which only plateaus at odd height
  have weight $y$ is
  \[
  \frac{(1-x)\left(A + \sqrt{A(A+4x^{2})}\right)}{2 x^{2} A},
  \]
  where $A = (1-x)(x^{2}(xy-x) +x-1)$.
\end{theorem}
\begin{proof}  Set $C_{2k+1} = xy - x$ and $C_{2k} =
  0$ in \eqref{eq:contfracexpansion2}. At even height, no correction is necessary. In this case, if we call the generating function $f(x,y)$, we
  easily get
  \[
  f(x, y) = \frac{1}{1-x-x^{2}(xy-x)-\cfrac{x^{2}}{1-x-x^{2} f(x,y)}};
  \]
  The proof is completed by solving for $f$.
\end{proof}

\emph{Remark}. The corresponding triangle is
  surprisingly close---but \emph{not} equal---to sequence A114581 in the
  OEIS. (They differ for paths of length $7$.)

\begin{theorem}
The generating function for Motzkin paths in which UHDs have weight $y$,
  UHHDs have weight $z$, and no plateaus of length three or more appear is
\begin{equation}
\label{pw}
  \frac{-(A(2) + 1) + \sqrt{(A(4)+1)(A(0) + 1)}}{2 x^2 (x-1)},
\end{equation}
  where $A(k) = x(x-1)(x^3z + x^2y +x + k)$.
\end{theorem}
\begin{proof}  In \eqref{eq:contfracexpansion2}, set $C_{k} = xy - x +
  x^{2}z - x^{2} - x^{3}/(1-x)$.  To explain why, notice that the $xy - x$ and $x^{2}z - x^{2}$ terms
  give the correct weights to UHD and UHHDs, respectively, and
  subtracting $x^{3}/(1-x) = x^{3} + x^{4} + x^{5} +\cdots$ eliminates
  the possibility of plateaus of length three or more. Using the
  functional equation~\eqref{newergfrelation} one can find that the
  generating function is as in \eqref{pw}.
  \end{proof}

\emph{Remark}. Notice that, by setting
  $y$ and $z$ to $1$ in \eqref{pw}, this is another way of counting the number of
  Motzkin paths with minimal plateau length 1 and maximal plateau length~2.
  (Compare with Prodinger and Wagner~\cite{prodinger}.)

\begin{theorem}
The generating function which counts Motzkin paths in which
  UHDs at height $2$ or more have weight~$y$, and UHHDs at a height that
  is a multiple of $3$ have weight~$z$ is given by
  \[
  G(x, y, z) = \frac{1}{1-x-x^{2}p},
  \]
where
\[
  p = -\frac{A+\sqrt{B}}{2x^{2}D},
  \]
  and
  \begin{align*}
    A &= \begin{gathered}[t]
         x^{2}Y(x^{2}Y(x^{2}Y+3x-3) + x^{2}Z(x^{2}Y + 2x-2) + 2x^{2} - 6x+3)\\
         {} + x^{2}Z(2x^{2} - 2x+1) - 2x^{2} + 3x - 1,
         \end{gathered}\\
    B &= \begin{gathered}[t]
         (x^{2}Y-1)(x^{2}Y+2x-1)
         (D - x^{3}Y - x^{3}Z - 2x^{2} +x)\\
         (D + x^{3}Z + x^{3}Y-x),
         \end{gathered}\\
    D &= x^{2}Y(x^{2}Z + x^{2}Y + x-2) + (y-z)x^{4} + x^{3}Z - 2x+1,\\
    Y &= xy - x, \\
    Z &= x^2z -x^2.
  \end{align*}
\end{theorem}
\begin{proof}  It is easy to see what our
  correction terms must be:

  \begin{center}
    \begin{tabular}{c|ccccccccccc}
      $k$    & 1 & 2   & 3     & 4  & 5 & 6 & 7 & 8 & 9 & 10 &  \\
      $C_{k}$ & 0 & $Y$ & $Y+Z$ & $Y$& $Y$ & $Y+Z$ & $Y$ & $Y$ & $Y+Z$ & $Y$& $\cdots$ \\
    \end{tabular}
  \end{center}

  Above, the correction term $Y$ equals $xy - x$ and $Z$ equals $x^2z -
  x^2$. If $G$ is the generating function for such paths, the continued
  fraction expansion for $g$ is
  \begin{multline}
  G = \hcfrac{1}{1-x} 
      - \hcfrac{x^{2}}{1-x - x^2 Y} 
      - \hcfrac{x^{2}}{1-x-x^{2} (Y + Z)} 
      - \hcfrac{x^{2}}{1-x-x^{2} Y}\\ 
      - \hcfrac{x^{2}}{1-x - x^2 Y} 
      - \hcfrac{x^{2}}{1-x-x^{2} (Y + Z)} 
      - \cdots
  \end{multline}
  We need to split off the ``purely periodic'' part of that continued
  fraction; call that generating function $p$. It satisfies
  \[
  p =  \hcfrac{1}{1-x-x^{2} Y} 
      - \hcfrac{x^{2}}{1-x-x^{2} (Y+Z)} 
      - \hcfrac{x^{2}}{1-x-x^2 Y - x^2 p} 
  \]
  and therefore equals
  \[
  p = -\frac{A+\sqrt{B}}{2x^{2}D},
  \]
  with $A$, $B$, and $D$ as in the statement of the theorem.  The theorem now follows.
  \end{proof}

\section{Acknowledgments}

We would like to thank JenAlyse Arena for getting us started on this
project. Also, many of the generating functions in this paper were
verified using the \path{sage.math.washington.edu} computer, which is
supported by National Science Foundation Grant No. DMS-0821725.

\renewcommand{\MR}[1]{\href{http://www.ams.org/mathscinet-getitem?mr=#1}{MR~#1}}
\providecommand{\ISBN}[1]{\href{http://worldcat.org/isbn/#1}{ISBN~#1}}
\bibliographystyle{amsplainurl}
\bibliography{references}

\end{document}